\documentclass[11pt,a4paper]{article}%
\usepackage{makeidx}
\usepackage{amssymb}
\usepackage{amsfonts}
\usepackage{graphicx}
\usepackage{amsmath}
\usepackage{fancyhdr}
\usepackage{geometry}
\usepackage{setspace}
\usepackage{color}
\usepackage{amsthm}%
\providecommand{\U}[1]{\protect\rule{.1in}{.1in}}
\newtheorem{thm}{Theorem}

\begin{document}

\title{Algebraic independence of certain infinite products involving\\the Fibonacci numbers
}
\author{\textsc{Daniel Duverney \thanks{ 
B{\^a}timent A1, 110
rue du chevalier Fran{\c{c}}ais, 59000 Lille, France \newline e-mail:
daniel.duverney@orange.fr},\,\,\,Yohei Tachiya \thanks{Hirosaki University,
Graduate School of Science and Technology, Hirosaki 036-8561, Japan \newline
e-mail: tachiya@hirosaki-u.ac.jp\newline This work was supported by JSPS
KAKENHI Grant Number JP18K03201} } }
\date{}
\maketitle

\begin{abstract}
Let $\{F_{n}\}_{n\geq0}$ be the sequence of the Fibonacci numbers. The aim of
this paper is to give explicit formulae for the infinite products
\[
\prod_{n=1}^{\infty}\left(  1+\frac{1}{F_{n}}\right)  ,\qquad\prod
_{n=3}^{\infty}\left(  1-\frac{1}{F_{n}}\right)
\]
in terms of the values of the Jacobi theta functions. From this we deduce the
algebraic independence over $\mathbb{Q}$ of the above numbers by applying
Bertrand's theorem
on the algebraic independence of the values of the Jacobi theta functions.



\end{abstract}

\thispagestyle{empty} \setcounter{page}{1}



\thispagestyle{empty} \setcounter{page}{1}



\pagestyle{fancy} \fancyhead{}
\fancyhead[EC]{Daniel Duverney and Yohei Tachiya} \fancyhead[EL,OR]{\thepage}
\fancyhead[OC]{
Algebraic independence of certain infinite products involving the Fibonacci numbers}
\fancyfoot{} \renewcommand\headrulewidth{0.5pt}


\section{Introduction and main results}

\label{sec1} Let $\{F_{n}\}_{n\geq0}$ be the Fibonacci sequence defined
recursively by
\[
F_{n+2}=F_{n+1}+F_{n},\qquad n\geq0
\]
with $F_{0}=F_{1}=1$. Arithmetical properties of infinite sums or products
involving the Fibonacci numbers have been investigated by several authors.
In~1989, Andr\'{e}-Jeannin~\cite{And} proved that the sum of the reciprocals
of the Fibonacci numbers $\sum_{n=1}^{\infty}1/F_{n}$ is irrational (see also
\cite{BV,Duv}). Moreover, it is shown in~\cite{DNNS} that the number
$\sum_{n=1}^{\infty}1/F_{2n-1}$ is transcendental. On the other hand, some
closed forms were discovered in particular cases; for example,
\begin{equation}
\sum_{n=1}^{\infty}\frac{1}{F_{2^{n}}}=\frac{5-\sqrt{5}}{2} \label{eq:3484}%
\end{equation}
(cf.~\cite[p.~225]{Luc}), which results from the use of telescoping series.
As for the infinite products,
the second author~\cite{Tac} derived that the numbers
\begin{equation}
\gamma_{j}:=\prod_{n=1}^{\infty}\left(  1+\frac{j}{F_{2^{n}}}\right)  ,\qquad
j=1,2,\dots\label{eq:397}%
\end{equation}
are all transcendental. The infinite products~\eqref{eq:397} can be regarded
as product analogues of \eqref{eq:3484}. The transcendence result on the
$\gamma$'s was extended in~\cite{KTT} to algebraic independence over
$\mathbb{Q}$ of the numbers $\gamma_{1}$, $\gamma_{2},\dots,\gamma_{m}$ for
any integer $m\geq1$. Also in the case of products, some closed forms can be
obtained through the telescoping method. For example, we have
\begin{equation}
\prod_{n=1}^{\infty}\left(  1+\frac{1}{F_{2n}}\right)  =1+\sqrt{5},\qquad
\prod_{n=2}^{\infty}\left(  1-\frac{1}{F_{2n}}\right)  =\frac{1+\sqrt{5}}{6},
\label{eq:453}%
\end{equation}
since
\begin{align*}
\prod_{n=1}^{\infty}\left(  1+\frac{1}{F_{2n}}\right)   &  =\prod
_{n=1}^{\infty}\frac{(1+\alpha^{-2(n-1)})(1-\alpha^{-2(n+1)})}{(1+\alpha
^{-2n})(1-\alpha^{-2n})}=\frac{2}{1-\alpha^{-2}}=2\alpha,\\[0.5em]
\prod_{n=2}^{\infty}\left(  1-\frac{1}{F_{2n}}\right)   &  =\prod
_{n=2}^{\infty}\frac{(1+\alpha^{-2(n+1)})(1-\alpha^{-2(n-1)})}{(1+\alpha
^{-2n})(1-\alpha^{-2n})}=\frac{1-\alpha^{-2}}{1+\alpha^{-4}}=\frac{\alpha}{3},
\end{align*}
respectively, where $\alpha:=(1+\sqrt{5})/2$ is the golden ratio. However, it
is still difficult to find closed forms and to determine arithmetical
properties for given infinite sums or products involving the Fibonacci
numbers. Note that we are still unaware of the transcendence of the numbers
$\sum_{n=1}^{\infty}1/F_{n}$ and $\sum_{n=1}^{\infty}1/F_{2n}$.

In this paper, we give explicit formulae for the two fundamental infinite
products
\[
\label{eq:90}\prod_{n=1}^{\infty}\left(  1+\frac{1}{F_{n}}\right)
=13.1509666577\dots,\qquad\prod_{n=3}^{\infty}\left(  1-\frac{1}{F_{n}%
}\right)  =0.1897891436\dots
\]
by means of the values of the Jacobi theta functions. Moreover, by using these
formulae we prove that the above two numbers are algebraically independent
over $\mathbb{Q}$. To state our results, we define the Jacobi theta functions
\begin{equation}
\vartheta_{2}(q):=2q^{1/4}\sum_{n=0}^{\infty}q^{n(n+1)},\qquad\vartheta
_{3}(q):=1+2\sum_{n=1}^{\infty}q^{n^{2}},\qquad\vartheta_{4}(q):=1+2\sum
_{n=1}^{\infty}(-1)^{n}q^{n^{2}}, \label{eq:0}%
\end{equation}
which converge for all complex number $q$ with $|q|<1$. Throughout this paper,
let
\begin{equation}
\beta:=\frac{1}{\alpha}=\frac{\sqrt{5}-1}{2}. \label{eq:3073}%
\end{equation}
Our main results are the following.

\begin{thm}
\label{thm:3} Let $\{F_{n}\}_{n\geq0}$ be the Fibonacci sequence. Then
\begin{align}
\xi_{1}  &  :=\prod_{n=1}^{\infty}\left(  1+\frac{1}{F_{n}}\right)
=2\beta^{-5/4}\frac{\vartheta_{2}(\beta)}{\vartheta_{4}(\beta^{4}%
),}\label{eq:12}\\[0.5em]
\xi_{2}  &  :=\prod_{n=3}^{\infty}\left(  1-\frac{1}{F_{n}}\right)
=\frac{\sqrt{5}}{6}\beta^{-5/4} \frac{\vartheta_{2}(\beta)\vartheta_{3}%
(\beta)\theta_{4}(\beta)}{\vartheta_{4}(\beta^{4})}. \label{eq:112}%
\end{align}

\end{thm}


\begin{thm}
\label{thm:33}The numbers $\xi_{1}$ and $\xi_{2}$ are algebraically
independent over $\mathbb{Q}$. In particular, the numbers $\xi_{1}$ and
$\xi_{2}$ are both transcendental.
\end{thm}


\section{Proofs of Theorems~\ref{thm:3} and \ref{thm:33}}

\label{Sec2} We first prove Theorem~\ref{thm:3}. Let $\vartheta_{2}%
(q),\vartheta_{3}(q),\vartheta_{4}(q)$ be the Jacobi theta functions defined
in \eqref{eq:0}. Using the triple-product identities
\[%
\begin{array}
[c]{c}%
\vartheta_{2}(q)=q^{1/4}\displaystyle\prod_{n=1}^{\infty}(1-q^{2n}%
)(1+q^{2n-2})(1+q^{2n}),\\
\vartheta_{3}(q)=\displaystyle\prod_{n=1}^{\infty}(1-q^{2n})(1+q^{2n-1}%
)^{2},\qquad\vartheta_{4}(q)=\displaystyle\prod_{n=1}^{\infty}(1-q^{2n}%
)(1-q^{2n-1})^{2}%
\end{array}
\]
(cf. \cite[Corollary~3.1]{BB}),
we have
\begin{align}
q^{-1/4}\vartheta_{2}(q)  &  =\prod_{n=1}^{\infty}(1-q^{2n})(1+q^{2n-2}%
)(1+q^{2n}),\label{eq:46}\\
\vartheta_{2}(q)\vartheta_{3}(q)\vartheta_{4}(q)  &  =2q^{1/4}\prod
_{n=1}^{\infty}(1-q^{2n})^{3},\label{eq:467}\\
\vartheta_{4}(q^{4})  &  =\prod_{n=1}^{\infty}(1-q^{4\cdot2n})(1-q^{4(2n-1)}%
)^{2}=\prod_{n=1}^{\infty}(1-q^{4n})(1-q^{4(2n-1)})\nonumber\\
&  =\prod_{n=1}^{\infty}(1-q^{2\cdot2n})(1-q^{2(2n-1)})(1+q^{2(2n-1)}%
)=\prod_{n=1}^{\infty}(1-q^{2n})(1+q^{4n-2}) \label{eq:47}%
\end{align}
(see~\cite[Proof of Proposition~3.1]{BB} for \eqref{eq:467}). Hence, we obtain
by \eqref{eq:46}, \eqref{eq:467}, and \eqref{eq:47}
\begin{align}
\prod_{n=1}^{\infty}\left(  1+\frac{1}{F_{2n-1}}\right)   &  =\prod
_{n=1}^{\infty}\frac{(1+\beta^{2n-2})(1+\beta^{2n})}{1+\beta^{4n-2}}%
=\beta^{-1/4}\frac{\vartheta_{2}(\beta)}{\vartheta_{4}(\beta^{4}%
)},\label{eq:09101}\\[0.5em]
\prod_{n=2}^{\infty}\left(  1-\frac{1}{F_{2n-1}}\right)   &  =\prod
_{n=2}^{\infty}\frac{(1-\beta^{2n-2})(1-\beta^{2n})}{1+\beta^{4n-2}}%
=\beta^{-1/4}\frac{1+\beta^{2}}{2(1-\beta^{2})}\frac{\vartheta_{2}%
(\beta)\vartheta_{3}(\beta)\vartheta_{4}(\beta)}{\vartheta_{4}(\beta^{4})}.
\label{eq:09102}%
\end{align}
Therefore, Theorem~\ref{thm:3} follows from \eqref{eq:3073}, \eqref{eq:09101},
\eqref{eq:09102}, and the formulae~\eqref{eq:453}.

Next we show Theorem~\ref{thm:33}. The well-known identities
\[
\vartheta_{3}^{4}(q)=\vartheta_{2}^{4}(q)+\vartheta_{4}^{4}(q),\qquad
2\vartheta_{3}^{2}(q^{2})=\vartheta_{3}^{2}(q)+\vartheta_{4}^{2}%
(q),\qquad\vartheta_{4}^{2}(q^{2})=\vartheta_{3}(q)\vartheta_{4}%
(q)\label{eq:2}%
\]
(cf. \cite[Chapter~2, \S ~2.1]{BB}) lead to
\begin{equation}
\vartheta_{2}^{4}(\beta)=\vartheta_{3}^{4}(\beta)-\vartheta_{4}^{4}%
(\beta),\qquad\vartheta_{4}^{4}(\beta^{4})=\frac{\vartheta_{3}(\beta
)\vartheta_{4}(\beta)}{2}(\vartheta_{3}^{2}(\beta)+\vartheta_{4}^{2}(\beta)).
\label{eq:492}%
\end{equation}
Hence, substituting \eqref{eq:492} into
\[
\xi_{1}^{4}=16\beta^{-5}\frac{\vartheta_{2}^{4}(\beta)}{\vartheta_{4}%
^{4}(\beta^{4})}%
\]
from \eqref{eq:12}, we have
\begin{equation}
\xi_{1}^{4}=32\beta^{-5}\frac{\vartheta_{3}^{2}(\beta)-\vartheta_{4}^{2}%
(\beta)}{\vartheta_{3}(\beta)\vartheta_{4}(\beta)}=32\beta^{-5}\left(
\frac{\vartheta_{3}(\beta)}{\vartheta_{4}(\beta)}-\frac{\vartheta_{4}(\beta
)}{\vartheta_{3}(\beta)}\right)  . \label{eq:09081}%
\end{equation}
On the other hand, by (\ref{eq:12}) and (\ref{eq:112}), we have also%
\begin{equation}
\vartheta_{3}(\beta)\vartheta_{4}(\beta)=\frac{12}{\sqrt{5}}\cdot\frac{\xi
_{2}}{\xi_{1}}. \label{eq:113}%
\end{equation}
By (\ref{eq:09081}) and (\ref{eq:113}) we see that the numbers $\vartheta
_{3}(\beta)/\vartheta_{4}(\beta)$ and $\vartheta_{3}(\beta)\vartheta_{4}%
(\beta)$ are algebraic over the field $\mathbb{K}:=\mathbb{Q}(\xi_{1},\xi
_{2})$, and thus so are the numbers $\vartheta_{3}(\beta)$ and $\vartheta
_{4}(\beta)$. Therefore, we can write
\[
2\geq\text{trans.deg}_{\mathbb{Q}}\mathbb{K}=\text{trans.deg}_{\mathbb{Q}%
}{\mathbb{K}}(\vartheta_{3}(\beta),\vartheta_{4}(\beta))\geq\text{trans.deg}%
_{\mathbb{Q}}\mathbb{Q}(\vartheta_{3}(\beta),\vartheta_{4}(\beta))=2,
\]
where at the last equality we used a result of Bertrand~\cite[Theorem~4]{Ber}
(see also \cite{EKT}) which is a consequence of Nesterenko's
theorem~\cite{Nes} on the algebraic independence of the values of Ramanujan
functions. This implies that $\text{trans.deg}_{\mathbb{Q}}\mathbb{K}=2$,
namely, the numbers $\xi_{1}$ and $\xi_{2}$ are algebraically independent
over~$\mathbb{Q}$. The proof of Theorem~\ref{thm:33} is completed.




\end{document}